\documentclass[a4paper]{amsart}

\usepackage{amssymb,latexsym}
\usepackage{amsmath}
\usepackage{amsfonts}
\usepackage{amsthm}
\usepackage{amssymb}

\theoremstyle{plain}
\newtheorem{theorem}{Theorem}
\newtheorem{lemma}[theorem]{Lemma}

\theoremstyle{definition}

\title[Rank zero curves induced by rational Diophantine triples]{Rank zero elliptic curves induced by rational Diophantine triples}
\begin{document}

\date{}


\author[A. Dujella]{Andrej Dujella}
\address[A. Dujella]{
Department of Mathematics\\
Faculty of Science\\
University of Zagreb\\
Bijeni{\v c}ka cesta 30, 10000 Zagreb, Croatia
}
\email{duje@math.hr}

\author[M. Miki\'c]{Miljen Miki\'c}
\address[M. Miki\'c]{Kumi\v{c}i\'ceva 20, 51000 Rijeka, Croatia
}
\email{miljen.mikic@gmail.com}

\begin{abstract}
Rational Diophantine triples, i.e. rationals $a,b,c$ with the property that $ab+1$, $ac+1$, $bc+1$
are perfect squares, are often used in construction of elliptic curves with high rank.
In this paper, we consider the opposite problem and ask how small can be the rank of
elliptic curves induced by rational Diophantine triples.
It is easy to find rational Diophantine triples with elements with mixed signs
which induce elliptic curves with rank $0$. However, the problem of finding
such examples of rational Diophantine triples with positive elements is much more challenging,
and we will provide the first such known example.
\end{abstract}


\maketitle

\section{Introduction}

A set $\{a_1, a_2, \dots, a_m\}$ of $m$ distinct nonzero rationals
is called {\it a rational Diophantine $m$-tuple} if $a_i a_j+1$ is a perfect square
for all $1\leq i < j\leq m$.
Diophantus discovered a rational Diophantine quadruple
$\{\frac{1}{16}, \frac{33}{16}, \frac{17}{4}, \frac{105}{16} \}$.
The first example of a Diophantine quadruple in integers, the set $\{1, 3, 8, 120\}$,
was found by Fermat. In 1969, Baker and Davenport \cite{B-D} proved that Fermat's set
cannot be extended to a Diophantine quintuple in integers.
Recently, He, Togb\'e and Ziegler proved that there are no Diophantine quintuples
in integers \cite{HTZ} (the nonexistence of Diophantine sextuples in integers was proved in \cite{D-crelle}).
Euler proved that there are infinitely many rational Diophantine quintuples.
The first example of a rational Diophantine sextuple, the
set $\{11/192, 35/192, 155/27, 512/27, 1235/48, 180873/16\}$, was found by Gibbs \cite{Gibbs1},
while Dujella, Kazalicki, Miki\'c and Szikszai \cite{DKMS} recently proved that there are infinitely
many rational Diophantine sextuples (see also \cite{Duje-Matija,DKP,DKP-split}). It is not known whether there
exists any rational Diophantine septuple.
For an overview of results on
Diophantine $m$-tuples and its generalizations see \cite{Duje-Notices}.

The problem of extendibility and existence of Diophantine $m$-tuples
is closely connected with the properties of the corresponding elliptic curves.
Let $\{a,b,c\}$ be a rational Diophantine triple. Then there exist nonnegative rationals
$r,s,t$ such that $ab+1=r^2$, $ac+1=s^2$ and $bc+1=t^2$.
In order to extend the triple $\{a,b,c\}$ to a quadruple,
we have to solve the system of equations
\begin{equation} \label{eq2}
ax+1= \square, \quad bx+1=\square, \quad cx+1=\square.
\end{equation}
We assign the following elliptic curve to the system \eqref{eq2}:
\begin{equation} \label{e1}
E:\qquad  y^2=(ax+1)(bx+1)(cx+1).
\end{equation}
We say that the elliptic curve $E$ is induced by the rational Diophantine triple $\{a,b,c\}$.

Since the curve $E$ contains three $2$-torsion points
$$ A=\Big( -\frac{1}{a}, 0 \Big), \quad B=\Big( -\frac{1}{b}, 0 \Big), \quad
C=\Big( -\frac{1}{c}, 0 \Big), $$
by Mazur's theorem \cite{Mazur}, there are at most four possibilities for the torsion group
over $\mathbb{Q}$ for such curves: $\mathbb{Z}/2\mathbb{Z} \times \mathbb{Z}/2\mathbb{Z}$,
$\mathbb{Z}/2\mathbb{Z} \times \mathbb{Z}/4\mathbb{Z}$,
$\mathbb{Z}/2\mathbb{Z} \times \mathbb{Z}/6\mathbb{Z}$
and $\mathbb{Z}/2\mathbb{Z} \times \mathbb{Z}/8\mathbb{Z}$.
In \cite{D-glasnik}, it was shown that all these torsion groups actually appear. Moreover, it was shown that every elliptic curve with torsion
group $\mathbb{Z}/2\mathbb{Z} \times \mathbb{Z}/8\mathbb{Z}$ is induced by a Diophantine triple
(see also \cite{CG}). Questions about the ranks of elliptic curves induced by Diophantine triples were
studied in several papers
(\cite{ADP, D-rocky, D-bordo,D-glasnik, D-JB-S, D-Peral-LMSJCM, D-Peral-RACSAM, D-Peral-JGEA,D-Peral-highrank}).
In particular, such curves were used
for finding elliptic curves with the largest known rank
over $\mathbb{Q}$ and $\mathbb{Q}(t)$ with torsion groups
$\mathbb{Z}/2\mathbb{Z} \times \mathbb{Z}/4\mathbb{Z}$
(\cite{D-Peral-LMSJCM,D-Peral-JGEA}) and
$\mathbb{Z}/2\mathbb{Z} \times \mathbb{Z}/6\mathbb{Z}$ (\cite{D-Peral-RACSAM}).

In this paper, we consider the question how small can be the rank of
elliptic curves induced by rational Diophantine triples.
We will see that it is easy to find rational Diophantine triples with elements with mixed signs
which induce elliptic curves with rank $0$ and that there exist such curves
with torsion groups $\mathbb{Z}/2\mathbb{Z} \times \mathbb{Z}/4\mathbb{Z}$,
$\mathbb{Z}/2\mathbb{Z} \times \mathbb{Z}/6\mathbb{Z}$ and
$\mathbb{Z}/2\mathbb{Z} \times \mathbb{Z}/8\mathbb{Z}$.
However, the problem of finding
such examples of rational Diophantine triples with positive elements is much harder, and they exist
only for torsion $\mathbb{Z}/2\mathbb{Z} \times \mathbb{Z}/8\mathbb{Z}$.
We describe the method for finding candidates for such curves
and by using {\tt magma} we are able to find one explicit example.

\section{Conditions for point $S$ to be of finite order}

Apart from three $2$-torsion points $A$, $B$ and $C$, the curve $E$
contains also the following two obvious rational points:
$$ P = (0, 1), \quad S = \Big( \frac{1}{abc}, \frac{rst}{abc} \Big). $$
It is not so obvious, but it is easy to verify that $S = 2R$, where
$$ R = \Big( \frac{rs + rt + st + 1}{abc}, \frac{(r + s)(r + t)(s + t)}{abc} \Big). $$
Thus, a necessary condition for $E$ to have the rank equal to $0$ is that
the points $P$ and $S$ have finite order.
The triple $\{a,b,c\}$ is regular, i.e. $c=a+b\pm 2r$ if and only if $S=\mp 2P$
(see \cite{D-bordo}).

By Mazur's theorem and the fact that $S\in 2E(\mathbb{Q})$, we have the following possibilities:
\begin{itemize}
\item $mP = \mathcal{O}$, $m=3,4,6,8$;
\item $mS = \mathcal{O}$, $m=2,3,4$.
\end{itemize}

In particular, since the point $P$ cannot be of order $2$, is it not possible to have simultaneously
rank equal to $0$ and torsion group $\mathbb{Z}/2\mathbb{Z} \times \mathbb{Z}/2\mathbb{Z}$.

By the coordinate transformation $x\mapsto \frac{x}{abc}$, $y\mapsto \frac{y}{abc}$, applied
to the curve $E$, we obtain the equivalent curve
\begin{equation} \label{e2}
E':\qquad  y^2=(x+ab)(x+ac)(x+bc),
\end{equation}
and the points $A$, $B$, $C$, $P$ and $S$ correspond to
$A'=(-bc,0)$, $B'=(-ac,0)$, $C'=(-ab,0)$, $P'=(0,abc)$ and $S'=(1,rst)$, respectively.
In the next lemma, we will investigate all possibilities for point $S$ to be of finite order.

\begin{lemma} \label{l:mS}
\mbox{}
\begin{itemize}
\item[(i)] The condition $2S=\mathcal{O}$ is equivalent to
$$ (ab+1)(ac+1)(bc+1)=0. $$
\item[(ii)] The condition $3S=\mathcal{O}$ is equivalent to
$$ 3+4(ab+ac+bc)+6abc(a+b+c)+12(abc)^{2}-(abc)^{2}(a^2 + b^2 + c^2 - 2ab - 2ac - 2bc)=0. $$
\item[(iii)] The point $S$ is of order $4$ if and only if
$$ ((ab+1)^2 - ab(c-a)(c-b))((ac+1)^2 - ab(c-a)(c-b))((bc+1)^2 - ab(c-a)(c-b)) = 0. $$
\end{itemize}
\end{lemma}

\proof
\mbox{}
\begin{itemize}
\item[(i)] The condition $2S'= \mathcal{O}$ implies $rst=-rst$, i.e. $rst=0$, and
$$ (ab+1)(ac+1)(bc+1)=0. $$

\item[(ii)]
From $3S'=\mathcal{O}$, i.e. $x(2S')=x(-S')=x(S')$,
the formulas for doubling points of elliptic curves give
\begin{align*}
&3+(ab+ac+bc)\\
&=\dfrac{9+4(ab+ac+bc)^{2}+(abc(a+b+c))^{2}+12(ab+ac+bc)}{4r^{2}s^{2}t^{2}}\\
&\ \ \ \ +\dfrac{6abc(a+b+c)+4abc(ab+ac+bc)(a+b+c)}{4r^{2}s^{2}t^{2}}.
\end{align*}
Thus we get
\begin{align*}
&4\left((abc)^{2}+abc(a+b+c)+(ab+ac+bc)+1\right)(3+ab+ac+bc)\\
&=9+12(ab+ac+bc)+\left(6abc(a+b+c)+4(ab+ac+bc)^{2}\right)\\
&\ \ \ \ \ \ +4abc(ab+ac+bc)(a+b+c)+\left(abc(a+b+c)\right)^{2},
\end{align*}
which is equivalent to
\begin{align*}
&3+4(ab+ac+bc)+6abc(a+b+c)+12(abc)^{2}\\
& \ \ \ \ \ \ \ \ - (abc)^{2}(a^{2} + b^{2} + c^{2} - 2ab - 2ac - 2bc)=0.
\end{align*}

\item[(iii)]
The condition that the point $S'$ is of order $4$ is equivalent to $2S'\in \{A',B',C'\}$.
Let us assume that $2S'=C'$ (other two cases are completely analogous).
From the formulas for doubling points of elliptic curves, we get
\begin{align*}
&2+(bc+ac)\\
&=\dfrac{9+4(ab+ac+bc)^{2}+\left(abc(a+b+c)\right)^{2}+12(ab+ac+bc)}{4r^{2}s^{2}t^{2}}\\
&\ \ \ \ +\dfrac{6abc(a+b+c)+4abc(ab+ac+bc)(a+b+c)}{4r^{2}s^{2}t^{2}},
\end{align*}
which is equivalent to
$$\left(1+2ab-abc(c-a-b)\right)^{2}=0, $$
or
$$(ab+1)^{2}=ab (c-a)(c-b) .$$
\end{itemize}
\qed

\section{Rank zero curves for triples with mixed signs}

Let us now consider three possibilities for $mS = \mathcal{O}$.

Assume first that $2S= \mathcal{O}$. By Lemma \ref{l:mS}(i), we have $(ab+1)(ac+1)(bc+1)=0$,
so we conclude that $a,b,c$ cannot have the same sign.
If we allow the mixed signs, then in this case we may assume that $b=-1/a$.
In \cite{D-glasnik}, the following parametrization of rational Diophantine triples of the form
$\{a,-1/a,c\}$ is given:
$$ a=\frac{ut+1}{t-u}, \quad b=\frac{u-t}{ut+1}, \quad c=\frac{4ut}{(ut+1)(t-u)}. $$
To find examples with rank $0$, let us assume that the triple $\{a,-1/a,c\}$ is regular.
This condition leads to $(u^2-1)(t^2-1)=0$, so we may take $u=1$.
If we take e.g. $t=2$, we obtain the curve with torsion group
$\mathbb{Z}/2\mathbb{Z} \times \mathbb{Z}/4\mathbb{Z}$ and rank $0$,
induced by the triple
$$ \Big\{3, -\frac{1}{3}, \frac{8}{3} \Big\}. $$

\medskip

Assume now that $3S= \mathcal{O}$.
If we also have $3P= \mathcal{O}$, then $P=\pm S$, a contradiction.
Hence, if the point $P$ has finite order, the only possibility that $P$ is of order $6$.
This implies $2P=\pm S$ and $c=a+b\mp 2r$.
By inserting $b=(r^2-1)/a$ and $c=a+b+2r$ in the condition from Lemma \ref{l:mS}(ii), we get
$$ (2ar-1+2r^2) (-a+2ar^2-2r+2r^3) (2a^2r-a-2r+4ar^2+2r^3)=0. $$
Thus,
$$ a= \frac{-2r(r^2-1)}{-1+2r^2}, \quad \mbox{or} \quad
\frac{-(-1+2r^2)}{2r}, \quad \mbox{or} \quad
\frac{1-4r^2 \pm \sqrt{1+8r^2}}{4r}. $$
Take
\begin{equation} \label{eq:z2z6a}
 (a,b,c)=\left(\frac{-2r(r-1)(r+1)}{-1+2r^2}, \frac{-(-1+2r^2)}{2r}, \frac{(-1+2r)(2r+1)}{2(-1+2r^2)r}\right).
\end{equation}
Then the condition $ab>0$ is equivalent to $r>1$ or $r<-1$, while the condition
$bc>0$ is equivalent to $-1/2 < r < 1/2$. Hence, $a,b,c$ cannot have the same sign.

The case
$$ (a,b,c)=\left(\frac{-(-1+2r^2)}{2r}, \frac{-2r(r-1)(r+1)}{-1+2r^2}, \frac{(-1+2r)(2r+1)}{2(-1+2r^2)r}\right) $$
is the same as the previous case, just $a$ and $b$ are exchanged.

Finally, let $8r^2+1 = (2rt+1)^2$, to get rid of a square root in the third case. It gives $r = \frac{-t}{-2+t^2}$. Then
$$ (a,b,c)=\left(\frac{-t(t-2)(t+2)}{2(-2+t^2)}, \frac{2(t-1)(t+1)}{(-2+t^2)t}, \frac{-(-2+t^2)}{2t}\right) $$
(or $a$ and $b$ exchanged).
The condition $ac>0$ is equivalent to $t>2$ or $t<-2$, while the condition
$bc>0$ is equivalent to  $-1 < t < 1$. Hence, is this case also $a,b,c$ cannot have the same sign.

If we allow the mixed signs, then we can obtain examples with rank $0$, e.g. from triples
of the form (\ref{eq:z2z6a}). E.g. for $t=4$ we obtain the curve with torsion group
$\mathbb{Z}/2\mathbb{Z} \times \mathbb{Z}/6\mathbb{Z}$ and rank $0$,
induced by the triple
$$ \Big\{-\frac{12}{7}, \frac{15}{28}, -\frac{7}{4} \Big\}. $$

\medskip

It remains the case when the point $S$ is of order $4$.
Then the point $R$, such that $2R=S$ is of order $8$ and therefore
the torsion group of $E$ is $\mathbb{Z}/2\mathbb{Z} \times \mathbb{Z}/8\mathbb{Z}$.
As we already mentioned in the introduction, it is shown in \cite{D-glasnik} that
every elliptic curve over $\mathbb{Q}$ with this torsion group
is induced by a rational Diophantine triple.
More precisely, any such curve is induced by a Diophantine triple of the form
\begin{equation}\label{z2z8T}
\left\{ \frac{2T}{T^2-1}, \,\, \frac{1-T^2}{2T}, \,\, \frac{6T^2-T^4-1}{2T(T^2-1)} \right\}.
\end{equation}
It is clear that the elements of (\ref{z2z8T}) have mixed signs.
By taking $T=2$ we obtain the curve with torsion group
$\mathbb{Z}/2\mathbb{Z} \times \mathbb{Z}/8\mathbb{Z}$ and rank $0$,
induced by the triple
$$ \Big\{ \frac{4}{3}, -\frac{3}{4}, \frac{7}{12} \Big\}. $$

\section{An example of rank zero curve for triple with positive elements}

In a previous section, we showed that for rational Diophantine triples with all positive elements
we cannot have rank $0$ and torsion group $\mathbb{Z}/2\mathbb{Z} \times \mathbb{Z}/4\mathbb{Z}$
or $\mathbb{Z}/2\mathbb{Z} \times \mathbb{Z}/6\mathbb{Z}$. So the only remaining possibility
is the torsion group $\mathbb{Z}/2\mathbb{Z} \times \mathbb{Z}/8\mathbb{Z}$.
Since the elements of (\ref{z2z8T}) clearly have mixed signs ($ab=-1$),
and all curves with torsion group $\mathbb{Z}/2\mathbb{Z} \times \mathbb{Z}/8\mathbb{Z}$
are induced by (\ref{z2z8T}), on the first sight we might think that
triples with all positive elements are not possible for this torsion group.
However, it is shown in \cite{D-Peral-RACSAM} that this is not true.
Namely, we may have a triple with positive elements which induce the same curve as
(\ref{z2z8T}) for certain rational number $T$.

But in \cite{D-Peral-RACSAM} it remained open whether it is possible to
obtain simultaneously torsion group $\mathbb{Z}/2\mathbb{Z} \times \mathbb{Z}/8\mathbb{Z}$
and rank $0$ for triples with positive elements, although some candidates for such triples
are mentioned.

As in the previous section, we assume that the point $S$ is of order $4$,
and we take $b=(r^2-1)/a$, $c=a+b+2r$. By inserting this in the first factor
$$ (ab+1)^2 - ab(c-a)(c-b) $$
in Lemma \ref{l:mS}(iii), we get the quadratic equation in $a$:
$$ (2r^3-2r)a^2+(4r^4-6r^2+1)a+2r^5+2r-4r^3 = 0.$$
Its discriminant,
$$ 1-4r^4+4r^2 $$
should be a perfect square.
The quartic curve defined by this equation is birationally equivalent to the elliptic curve
$$ E_1: \quad Y^2 = X^3+X^2+X+1 $$
with rank $1$ and a generator $P_1=(0, 1)$.
Thus, by computing multiples of the point $P_1$ on the curve $E_1$
(adding the $2$-torsion point $T_1=(-1,0)$ has a same effect as changing $r$ to $-r$),
and transferring them back to the quartic, we obtain candidates for the solution of our problem.
However, we have to satisfy the condition that all elements of the corresponding triple are positive
(it is enough that all elements have the same sign, since by multiplying all elements from
a rational Diophantine triple by $-1$ we obtain again a rational Diophantine triple).
The first two multiples of $P$ producing the
triples with positive elements are $6P$ and $11P$.

The point $6P$ gives
$r=-\frac{3855558}{3603685}$ and the triple
$$ (a,b,c) = \left(
\frac{1884586446094351}{25415891646864180},
\frac{14442883687791636}{7402559392524605},
\frac{60340495895762708555}{14487505263205637124} \right). $$
We were not able to determine that the rank of the correspoding curve.
Namely, both {\tt magma} and {\tt mwrank}
give that $0 \leq \mbox{rank} \leq 2$. Assuming the Parity conjecture the rank should be equal to $0$ or $2$.

The point $11P$ gives
$r=\frac{35569516882766685106979}{32383819387240952672281}$ and the triple $(a,b,c)$, where
\begin{align*}
a&= \frac{69705492951192675600645567228019184577147632882703132983}{132014843349912467692901303836561266921302184459536763120}, \\ b&= \frac{47826829880079829075801189563942620732062701095548790400}{122336669420709509303637442647966391336596694969835459327}, \\ c&=
 \frac{47982111146649404421749331709393501777791774558546217987550257759801}{15400090753918257364093484910580652390786084055043677020804056653840}.
\end{align*}
(By comparing $j$-invariants, we get that the same curve is induced by (\ref{z2z8T}) for
$T=\frac{18451786408106133183649}{41916048174422594852689}$.)
For the corresponding curve, both {\tt mwrank} and {\tt magma} function {\tt MordellWeilShaInformation}
give that $0 \leq \mbox{rank} \leq 4$.
However, {\tt magma} (version  {\tt V2.24-7}) function {\tt TwoPowerIsogenyDescentRankBound},
which implements the algorithm by Fisher from \cite{Fisher}, gives that the rank is equal to $0$
(at step $5$, just beyond $4$-descent, but not yet $8$-descent).
Hence, we found an example of a rational Diophantine triple with positive elements for
which the induced elliptic curve has the rank equal to $0$. Let us mention that the same {\tt magma} function applied to the
curve mentioned above corresponding to the point $6P$ gives only $\mbox{rank} \leq 2$.
This construction certainly gives
infinitely many multiples of $P$ which produce triples with positive elements
(the set $E_1(\mathbb{Q})$ is dense in $E_1(\mathbb{R})$, see e.g. \cite[p.78]{Skolem}).
However, it is hard to predict distribution of ranks in such families of elliptic curve,
so we may just speculate that there might be infinitely many curves in this family with rank $0$.

\bigskip

{\bf Acknowledgements.}
A.D. was supported by the Croatian Science Foundation under the project no.~IP-2018-01-1313.
He also acknowledges support from the QuantiXLie Center of Excellence, a project
co-financed by the Croatian Government and European Union through the
European Regional Development Fund - the Competitiveness and Cohesion
Operational Programme (Grant KK.01.1.1.01.0004).

\end{document}